

\input amstex
\documentstyle{amsppt}
\magnification=1200
\catcode`\@=11
\redefine\logo@{}
\catcode`\@=13

\define \bn{\Bbb N}
\define \bz{\Bbb Z}
\define \bq{\Bbb Q}
\define \br{\Bbb R}
\define \bc{\Bbb C}

\define \M{{\Cal M}}
\define\Ha{{\Cal H}}
\define\La{{\Cal L}}  
 
\define\0o{{\overline 0}}
\define\1o{{\overline 1}}
\define\rk{\text{rk}~}


\define\mult{\text{mult}}

\define\pd#1#2{\dfrac{\partial#1}{\partial#2}}

\define\NEF{\text{NEF}}
\define\Aut{\text{Aut}}

\TagsOnRight

\document

\topmatter
\title
The Arithmetic Mirror Symmetry and Calabi--Yau manifolds 
\endtitle

\author
Valeri A. Gritsenko \footnote{Supported by
RIMS of Kyoto University 
\hfill\hfill}
and
Viacheslav V. Nikulin \footnote{Supported by
Grant of Russian Fund of Fundamental Research and RIMS of Kyoto University 
\hfill\hfill}
\endauthor

\address
St. Petersburg Department of Steklov Mathematical Institute,
\newline
${}\hskip 8pt $
Fontanka 27, 191011 St. Petersburg,  Russia
\endaddress
\email
gritsenk\@gauss.pdmi.ras.ru 
\endemail

\address
Steklov Mathematical Institute,
ul. Vavilova 42, Moscow 117966, GSP-1, Russia
\endaddress
\email
slava\@nikulin.mian.su
\endemail

\abstract
We extend our variant of mirror symmetry for K3 
surfaces \cite{GN3} and clarify its relation with mirror symmetry for 
Calabi--Yau manifolds. We introduce two classes (for the models A and B) of 
Calabi--Yau manifolds fibrated by K3 surfaces with some special 
Picard lattices. These two classes are 
related with automorphic forms on IV type domains which we  
studied in our papers \cite{GN1}---\cite{GN6}. Conjecturally these 
automorphic forms take part in the 
quantum intersection pairing for model A, 
Yukawa coupling for model B and mirror symmetry between these two 
classes of Calabi--Yau manifolds. 
Recently there were several papers by physicists where it was shown 
on  some examples. We propose   
a problem of classification of introduced 
Calabi--Yau manifolds.  
Our papers \cite{GN1}---\cite{GN6} and \cite{N3}---\cite{N14} 
give a hope that this is possible. They describe possible 
Picard or transcendental lattices of general  
K3 fibers of the Calabi--Yau manifolds.   
\endabstract

\rightheadtext
{Arithmetic Mirror Symmetry}
\leftheadtext{V. Gritsenko and  V. Nikulin}
\endtopmatter

\head
0. Introduction
\endhead

In \cite{GN3} we suggested a variant of mirror symmetry for K3 surfaces 
which is related with reflection groups in hyperbolic spaces and 
automorphic forms on IV type domains. This subject was developed 
in our papers \cite{GN1}---\cite{GN6}, \cite{N11}---\cite{N14}.  
Some results of R. Borcherds \cite{B1}---\cite{B7} are also connected  
with this subject.  

Recently several papers by physicists have appeared where 
our automorphic forms (and some automorphic forms 
constructed by R. Borcherds) were used in mirror symmetry 
for Calabi--Yau manifolds. Physicists have shown that automorphic forms on 
IV type domains which we considered for our variant \cite{GN3} 
of mirror symmetry for K3 
surfaces take part in the 
quantum intersection pairing and the Yukawa coupling for 
some Calabi--Yau manifolds. 
We only mention papers which are directly  connected 
with this subject: Harvey -- Moore \cite{HM1}---\cite{HM3};  
Henningson -- Moore \cite{HeM1}, \cite{HeM2}; 
Kawai \cite{Ka1}, \cite{Ka2}; 
Dijkgraaf -- E. Verlinde -- H. Verlinde \cite{DVV}; 
Cardoso -- Curio -- L\"ust \cite{CCL}, 
but there are many other papers which are connected indirectly. 

For further study it is important to clarify  
relation between our variant of mirror symmetry for K3 surfaces,  
and the Calabi--Yau manifolds which 
were used by physicists. We want to give a definition of 
Calabi--Yau manifolds related with 
the variant \cite{GN3} of mirror symmetry for K3 surfaces  
and propose a problem of their classification. 
Our papers \cite{GN1}---\cite{GN6} and \cite{N3}---\cite{N14} 
give a hope that this is possible. They describe possible Picard 
or transcendental lattices of general K3 fibers of 
these Calabi--Yau manifolds. 

This paper was subject of the talk given by 
the second author at the conference ``Mirror Symmetry and 
Calabi--Yau manifolds'' held at RIMS Kyoto University  
on 2---6 December, 1996. It  was written during our stay 
at RIMS on 1996---1997.  
We are grateful to the Mathematical Institute for hospitality. 

The preliminary variant of this paper was published as preprint 
\cite{GN7}. 

\head
1. K3 surfaces with skeleton as fibers of Calabi--Yau manifolds 
for model A 
\endhead 

The main property of 
Calabi--Yau manifolds $\Cal X$ which are related with our variant 
of mirror symmetry for K3 surfaces 
and which were considered in papers of 
physicists we cited above is that {\it they are fibrated by K3 
surfaces}. There is a morphism 
$$
\pi :{\Cal X} \to B
$$
with the general fiber which is a K3 surface $X$. Let $S$ be 
the Picard lattice of the general fiber of $\pi$. Then 
${\Cal X}$ is related with the moduli $\M_S$ of K3 surfaces with 
the Picard lattice $S$. It is known that $\M_S=G\setminus \Omega(T)$ where 
$\Omega(T)$ is a Hermitian symmetric domain of type IV which is 
defined by the transcendental lattice $T=S^\perp_{L_{K3}}$. Here 
$L_{K3}$ is the second cohomology lattice of K3 and 
$G\subset O(T)$ is a subgroup of finite index. Thus, K3 fibrated 
Calabi--Yau manifolds are related with arithmetic quotients  
$G\setminus \Omega (T)$ of IV type domains. 
One can pick up the family $\pi$ above from a morphism  
$B \to \overline{\M_S}=\overline{G \setminus \Omega (T)}$.  

Recently there was a big progress in studying mirror symmetry for 
Calabi--Yau complete intersections in toric varieties. We mention papers of   
Candelas -- de la Ossa -- Green -- Parkers \cite{COGP}, Morrison 
\cite{M}, Batyrev \cite{Bat}, Kontsevich \cite{Ko}, 
Kontsevich -- Manin 
\cite{KoM}, Ruan -- Tian \cite{RT} and Givental \cite{Gi}. 
For these Calabi--Yau manifolds  
mirror symmetry comes from duality between polyhedra defining ambient 
toric varieties. Moreover, the Yukawa coupling and the 
quantum intersection pairing of Calabi--Yau compete intersections 
in toric varieties are strongly related with quantum cohomology 
of ambient toric varieties.    
One can ask about similar theory for K3 fibrated 
Calabi--Yau manifolds when one replaces toric varieties by arithmetic 
quotients of IV type domains defined by moduli of K3 surfaces with 
condition on Picard lattice. 
One can expect that 
{\it in some cases mirror symmetry for K3 fibrated 
Calabi--Yau manifolds is dominated by some variant of mirror symmetry 
for their K3 fibers}. This variant has been 
suggested in \cite{GN3}. Below we only extend it a little. 

\definition{Definition 1} An algebraic  
K3 surface $X$ over $\bc$ has {\it skeleton} if  
there exists an integral non-zero nef element $r \in \NEF (X)\cap S$ 
such that $r$ is invariant with respect to $\Aut (X)$. 
The element $r$ is called 
{\it canonical nef element}. Here we denote by  
$\NEF(X)\subset S\otimes \br$ the {\it nef cone}  
(equivalently, the clouser of the {\it K\"ahler cone}) 
of $X$ and by $S$ the {\it Picard lattice} of $X$.  
\enddefinition

Let $W^{(2)}(S)\subset O(S)$ be the group generated by reflections 
in all elements with square $-2$ of $S$ and $P(X)$ the 
set of all irreducible non-singular rational curves on $X$. 
From the global Torelli Theorem for K3 surfaces and 
description of the automorphism groups of K3 surfaces 
\cite{P-\u S\u S}, it follows that 
the canonical map $\pi: \Aut(X)\to \Aut(\NEF(X))$ 
has finite kernel and cokernel. Here $\Aut(\NEF(X))=\{\phi \in 
O(S)\,|\,\phi(\NEF(X))=\NEF(X)\}$.  
Moreover, the $\M=\NEF(X)/\br_{++}$ is 
a fundamental chamber for the group $W^{(2)}(S)$ acting in 
the the hyperbolic space $\La(S)$ defined by $S$ with the 
set $P(X)$ of orthogonal vectors to faces (of highest dimension) 
of $\NEF(X)/\br_{++}$. 
We then get the following description of K3 surfaces $X$ with skeleton 
and their Picard lattices $S$.  

a) {\it Elliptic type.} The set $P(X)$ is finite 
and generates $S\otimes \bq$. Then one can 
find a canonical nef element $r$ with $r^2>0$. One can consider the finite 
set of non-singular rational curves on $X$ as a ``skeleton'' of $X$. It 
essentially defines geometry of $X$. A K3 surface $X$ has elliptic 
type if and only if the group $W^{(2)}(S)$ has finite index in $O(S)$, 
and for $\rk S=2$ the $S$ contains at least $4$ different elements 
with square $-2$.      
Then $S$ is called {\it $2$-reflective of elliptic type}. 

a$^\prime$) {\it Special elliptic type.}  
The Picard lattice $S$ is one-dimensional. Then $r^2>0$. 
The lattice $S$ is called {\it $2$-reflective of special elliptic type}.  

b) {\it Parabolic type.} There exists a {\it canonical} elliptic fibration 
$|r|:X \to {\Bbb P}^1$ which is preserved by a subgroup 
$G\subset \Aut (X)$ of finite index. If $\Aut (X)$ is infinite, then 
the canonical elliptic fibration $|r|$ is unique and is 
preserved by the $\Aut (X)$ which is a finitely 
generated Abelian group up to finite index.     
For this case non-singular rational curves $P(X)$ on $X$ have 
bounded degree with respect to $r$. The $r$ together with $P(X)$ 
also can be considered as a ``skeleton'' of $X$. 
A K3 surface $X$ has parabolic type if and 
only if the quotient group $O^+(S)/W^{(2)}(S)$ considered as a group 
$\Aut (\M)$ of symmetries of a fundamental polyhedron $\M$ for 
$W^{(2)}(S)$ contains a subgroup $G\subset \Aut (\M)$ of 
finite index which fixes a non-zero  
element $r \in S$ with $r^2=0$ (it is unique and is 
fixed by $\Aut (\M)$ if $\Aut (\M )$ is infinite).  
Then $S$ is called {\it $2$-reflective of parabolic type}.   

\smallpagebreak 

Hyperbolic lattices $S$ satisfying one of the 
conditions above are called {\it 2-reflective}. 
Thus, {\it a K3 surface $X$ has skeleton if and only if 
its Picard lattice is $2$-reflective (of elliptic, including 
special elliptic, or parabolic type).}    
All K3 surfaces $X$ with skeleton are distributed 
in several families $\M_S$ of K3 surfaces with condition on Picard 
lattice according to their Picard lattice $S$. Any $S$ of rank one is 
$2$-reflective. If $\rk S=2$, the lattice $S$ is $2$-reflective if 
and only if it has a non-zero element with norm $-2$ or $0$.  
It was proved in 
\cite{N3}, \cite{N4}, \cite{N12} that for $\rk S\ge 3$ the set of 
2-reflective hyperbolic lattices $S$ is finite. All 
2-reflective hyperbolic lattices $S$ of elliptic type 
were classified in \cite{N3}, \cite{N7}, \cite{N8}. 

In practice, to construct K3 surfaces $X$ with skeleton, one need to 
find a polyhedron 
(the fundamental polyhedron 
$\M=\NEF(X)/\br_{++}$ for $W^{(2)}(S)$) in hyperbolic 
space with some condition of finiteness of volume and some condition 
of integrity (its Gram matrix should define 
a symmetric generalized Cartan matrix (e.g. see \cite{K})).  
Thus, the theory of K3 surfaces 
with skeleton is in fact 
{\it similar to the theory of toric varieties} where 
one need to consider some polyhedra in Euclidean space with a lattice.  

{\it For the model A of mirror symmetry we consider Calabi--Yau manifolds 
fibrated by K3 surfaces with skeleton (equivalently, with a 
2-reflective Picard lattice). Here and in what follows by the 
K3 surface we always mean a general K3 fiber of a  
general Calabi--Yau manifold.} From this point of view, 
it is a very interesting problem to classify Calabi--Yau 
manifolds fibrated by K3 surfaces with skeleton.   

\proclaim{Problem A}  Find all 2-reflective hyperbolic 
lattices $S$ such that there exists a Calabi--Yau manifold of dimension 
$\ge 3$ fibrated by K3 surfaces $X$ with the Picard lattice $S$. Find for 
this $S$ all Calabi--Yau manifolds fibrated by K3 surfaces with the 
Picard lattice $S$.  
\endproclaim

This problem looks much easier than the general problem of 
classification of Calabi--Yau manifolds fibrated by K3 surfaces because 
Picard lattices $S$ of K3 surfaces with skeleton are very special 
(and actually known). The monodromy group on the Picard 
lattice is easy to control. These Calabi--Yau manifolds 
have very special divisors defined by non-singular 
rational curves in K3 fibers.    

There are some examples of Calabi--Yau manifolds fibrated by  
K3 surfaces $X$ with skeleton. 
For example hyperbolic lattices $S$ with a 2-elementary discriminant 
group: $S^\ast /S\cong (\bz/2\bz)^a$, give a big part of 
$2$-reflective hyperbolic lattices $S$ of high rank and correspond 
to K3 surfaces with non-symplectic involutions (see \cite{N3} and 
also \cite{N8}). 
Using these involutions, C. Borcea \cite{Bo} and Cl. Voisin \cite{V} 
constructed Calabi--Yau $3$-folds fibrated by K3 surfaces with 
these Picard lattices. One can construct some other  
examples as complete intersections in toric varieties. 
But even for $\rk S \ge 3$ not for all 2-reflective 
lattices $S$ we know existence of Calabi--Yau manifolds 
fibrated by K3 surfaces with the Picard lattice $S$.

We remark that the theory of reflective hyperbolic lattices of elliptic 
and parabolic type was extended to {\it hyperbolic type} 
in \cite{N14}. One can introduce K3 surfaces with skeleton 
(or with $2$-reflective Picard lattice) of hyperbolic type. 
Without any doubt they are also important for mirror symmetry.   

\head
2. $2$-reflective automorphic forms, and transcendental lattices 
for K3 fibers of Calabi--Yau manifolds of model B
\endhead 

Now we consider the model B for mirror symmetry.  

The mirror symmetric subject to 2-reflective 
hyperbolic lattices $S$ is given by so called 
$2$-reflective lattices $T$ with 2 positive 
squares and by so called $2$-reflective automorphic forms on 
the symmetric domain 
$$
\Omega (T)=\{\bc \omega \subset T\otimes \bc \ |\ 
\omega \cdot \omega=0, \ \omega \cdot \overline{\omega} > 0 \}_0 .   
$$
Here a holomorphic automorphic form $\Phi$ on $\Omega (T)$  
with respect to a subgroup of $O(T)$ of finite index is called 
{\it $2$-reflective for $T$} if the divisor of $\Phi$ in 
$\Omega (T)$ is union of 
quadratic divisors $\Ha_\delta =\{\bc \omega \in \Omega (T)\ |\ 
\omega \cdot \delta=0\}$ (with some multiplicities)  
orthogonal to elements $\delta \in T$ with $\delta^2=-2$. A 
lattice $T$ having a $2$-reflective automorphic form is called 
{\it $2$-reflective}.   

Geometrical meaning of a $2$-reflective lattice $T$ with two 
positive squares  
and a $2$-reflective 
automorphic form $\Phi$ of $T$  
is that $\Phi$ is  
equal to zero only on the discriminant of the 
moduli $\Cal M_{T^\perp}$ of K3 surfaces with the 
Picard lattice $T^\perp$ (or with the transcendental lattice $T$).  
Here we use the global Torelli theorem \cite{P-\u SS} and epimorphicity 
of the Torelli map \cite{Ku}. 
All $2$-reflective automorphic 
forms $\Phi$ corresponding to $T$ define a semi-group which is 
very interesting. These automorphic forms take 
part in the mirror symmetry which we will consider. 

We expect that like the set of 
$2$-reflective hyperbolic lattices $S$ the set 
of $2$-reflective lattices $T$ with two positive squares 
is very small. Here the main conjecture is 
(see \cite{N13} and \cite{GN5} for more general 
and exact formulation) 

\proclaim{Arithmetic Mirror Symmetry Conjecture } 

a) 
The set of $2$-reflective 
lattices $T$ of $\rk T\ge 5$ is finite.

b) For any primitive isotropic $c \in T$ such that $c^\perp_T$ 
contains an element with square $-2$, the hyperbolic lattice 
$S=c^\perp_T/\bz c$ is 2-reflective. 
\endproclaim 

We will explain why we suppose that this is true. For a lattice $T$ 
we denote by $\Delta^{(2)}(T)$ 
the set of all elements of $T$ with square $-2$ and 
by $\Ha_\delta\subset \Omega (T)$ a quadratic divisor 
$\Ha_\delta$ which is orthogonal to $\delta \in \Delta^{(2)}(T)$. 
By Koecher principle (e. g. see \cite{Ba}),  
any automorphic form on a IV type domain $\Omega (T)$ has zeros 
if $\text{codim}_{\Omega (T)}{\Omega (T)_\infty}\ge 2$. Considering 
restriction of a $2$-reflective automorphic form $\Phi$ to 
subdomains $\Omega (T_1)$ where $T_1\subset T$, we get that 
$$
\left( \bigcup_{\delta \in \Delta^{(2)} (T)}{\Ha_\delta}\right) 
\bigcap \Omega (T_1)\not=
\emptyset
\tag{1}
$$
for any primitive sublattice $T_1\subset T$ with two positive 
squares such that $\text{codim}_{\Omega (T_1)}{\Omega (T_1)_\infty}$ $\ge 2$. 
One can consider condition \thetag{1} as the {\it analog of condition 
of finiteness of volume} for  
a polyhedron in hyperbolic space. This condition 
is extremely strong (see \cite{N13}), 
and we expect that lattices $T$ satisfying \thetag{1} 
satisfy the arithmetic mirror symmetry conjecture. Using \thetag{1}, 
it was shown in \cite{N13}   

\proclaim{Theorem 1} Let $T_n$ be the transcendental lattice of a general 
(i.e. with the Picard number 1) algebraic K3 surface of degree $n$. 
For any $N>0$ there exists $n>N$ 
such that the lattice $T_n$ is not $2$-reflective. In particular, 
the discriminant of moduli $\M_n$ of general K3 surfaces of the degree $n$ 
is not equal to zero set of any automorphic form on $\Omega (T_n)$.  
\endproclaim 

Arithmetic mirror symmetry Conjecture is very important for classification 
of $2$-reflec\-tive lattices $T$. After classification of 
$2$-reflective hyperbolic lattices $S$ (their set is finite for 
$\rk S \ge 3$), using arithmetic mirror symmetry Conjecture 
it would not be difficult to find all $2$-reflective lattices 
$T$ with $2$ positive squares. 

{\it For the model B of mirror symmetry we consider Calabi--Yau manifolds 
fibrated by K3 surfaces with a $2$-reflective transcendental lattice $T$. 
Here and in what follows by the transcendental lattice 
we always mean the transcendental lattice of 
a general K3 fiber of a general Calabi--Yau manifold}.   
Similarly to Problem A we propose 

\proclaim{Problem B} Find all $2$-reflective lattices $T$ 
such that there exists a Calabi--Yau manifold fibrated 
by K3 surfaces with the transcendental lattice $T$. Classify all 
Calabi--Yau manifolds fibrated by K3 surfaces 
with the $2$-reflective transcendental lattice $T$.
\endproclaim

Like for Problem A, we expect that classification of Calabi--Yau manifolds 
fibrated by K3 surfaces 
with a $2$-reflective transcendental lattice $T$ is much 
simpler than the general problem of classification of K3 fibrated 
Calabi--Yau manifolds because the set of $2$-reflective transcendental 
lattices $T$ is very small.  

\head
3. Mirror Symmetry
\endhead 

Let $T$ be a $2$-reflective lattice and $\Phi$ a $2$-reflective 
automorphic form of $T$. Considering product of $g^\ast \Phi$ over all  
$g \in G\setminus O^+(T)$, we can suppose that $\Phi$ is automorphic 
with respect to $O^+(T)$. For this case it is expected that $\Phi$ 
has a very special Fourier expansion at cusps $c \in T$.  
Here $c\in T$ is a primitive non-zero element with $c^2=0$.  

For simplicity we suppose that $T=U(k)\oplus S$ and $c \in U(k)$ 
(general case can be treated like in \cite{GN5, Sect. 2.3}).  
Here $\{ c, e\}$ is a basis of $U(k)$ such that $c^2=e^2=0$ and 
$c\cdot e=k \in \bn$.  
We consider the {\it mirror symmetry coordinate}  
$z \in \Omega (V^+(S))=S\otimes \br + i V^+(S)$ where $V^+(S)$ is 
the light cone of $S$. We associate to $z\in V^+(S)$ the point 
$$
\bc \omega_0 \in \Omega (T),\ \ \ \omega_0=((-z^2/2)c+(1/k)e)\oplus z.
$$  
Here $\omega_0 \in \bc \omega_0\in \Omega (T)$ is chosen by the 
condition $\omega_0\cdot c=1$. It is the mirror symmetry normalization. 
We expect that after identification of $S$ with the Picard 
lattice of a K3 surface $X$ the {\it $2$-reflective 
automorphic form $\Phi$} 
multiplied by  some constant {\it could be written as}   
$$
\split
&\Phi (z)=\\
&=\sum_{w \in W^{(2)}(S)}{\epsilon (w)
\left(\exp{(2\pi i (w(\rho)\cdot z))}-
\sum_{a \in \text{NEF}(S)}{N(a)\exp{(2\pi i (w(\rho+a)\cdot z))}}\right)}\\
&=\exp{\left(2\pi i (\rho\cdot z)\right)}
\prod_{\alpha \in \text{EF}(S)}{\left(1-
\exp{(2\pi i(\alpha\cdot z))}
\right)^{\mult~\alpha}}. 
\endsplit
\tag{2}
$$   
Here $\epsilon:W^{(2)}(S)\to \{\pm 1\}$ is some character 
and $\text{EF}(S)$ is the set of effective elements of $X$.  
All Fourier coefficients $N(a)$ and 
``multiplicities'' $\mult~\alpha$ 
are integral. The $\rho$ is a non-zero element of $\bq_+ \cdot \NEF(S)$. 
It is called the {\it generalized lattice Weyl vector}. 
The infinite 
product in this formula 
is the product of Borcherds type \cite{B5}.   
We remark that {\it if one has the Fourier expansion of type 
\thetag{2} with a non-trivial reflection group $W^{(2)}(S)$, then 
the lattice $S$ is automatically $2$-reflective} (see \cite{N12}, 
\cite{GN5}). It is why we were previously 
forced to restrict by $2$-reflective 
hyperbolic lattices $S$ for the model A. Existence of the 
Fourier and the infinite product expansion of type \thetag{2} 
is also important in  Physics. E. g. see \cite{CCL}, 
\cite{DVV}, \cite{HM1}---\cite{HM3}, \cite{HeM1}, \cite{HeM2}, 
\cite{Ka1}, \cite{Ka2}.   
  
We consider the families $\M_{T^\perp}$ and $\M_S$ of K3 surfaces 
with condition on Picard lattice as 
{\it mirror symmetric} if for $T$ and $S=c^\perp_T/\bz c$  
a $2$-reflective automorphic form $\Phi$ satisfying \thetag{2} 
does exist. Remark that 
here both lattices $T$ and $S$ (if the group $W^{(2)}(S)$ 
is non-trivial) are $2$-reflective. This definition 
of mirror symmetry for K3 surfaces 
was used in \cite{GN3} for some more narrow class of $\Phi$.  

Geometrically existence of the form $\Phi$ 
is very nice. On the one hand, $\Phi$ 
is an ``algebraic function" on the moduli $\M_{T^\perp}$ 
with zeros only on the discriminant, and the identity 
\thetag{2} reflects geometry of moduli $\M_{T^\perp}$ 
in the neighborhood of the cusp $c$. On the other hand, 
the identity (2) reflects geometry of curves on 
general K3 surfaces $X \in \M_S$ of 
the mirror symmetric family. 
It is why we considered in \cite{GN3} this 
type of mirror symmetry for K3 surfaces 
as a very natural and beautiful one. 
It seems that the case when $\Phi$ has 
zeros of multiplicity one is the most important. 
Then one can associate to \thetag{2} the   
{\it generalized Lorentzian Kac--Moody superalgebra with 
the denominator function \thetag{2}} 
(this case was considered in \cite{GN3}).  
Moreover, it seems important to have $\Phi$ which 
is equal to zero along all quadratic divisors orthogonal 
to elements of $T$ with square $-2$ and with multiplicities 
which are as small as possible.  

Like for toric geometry, we can introduce    
{\it K3 surfaces $X$ with reflexive 2-reflective Picard 
(or transcendental) lattice}   
when $X$ takes part in our mirror symmetry on the (A) or (B) 
side respectively, or both sides.     
Finding of all these cases is especially interesting. 

\proclaim{Problems A$^\prime$, B$^\prime$ } 
Find all reflexive 2-reflective lattices $S$ (respectively $T$) 
such that there exists a Calabi--Yau manifold 
fibrated by K3 surfaces with the Picard lattice $S$ (respectively 
with the transcendental lattice $T$). Classify all Calabi--Yau manifolds 
fibrated by K3 surfaces with the reflexive 2-reflective Picard 
lattice $S$ (respectively with the reflexive 2-reflective transcendental 
lattice $T$). 
\endproclaim 

The first multi-dimensional automorphic form of type \thetag{2} 
was found by R. Bor\-cherds \cite{B2}. 
He constructed it for   
the even unimodular lattice $T$ of rank $28$. Then $S$ 
is the even unimodular hyperbolic lattice of rank $26$. It is 
$2$-reflective of parabolic type.  
It is expected that this 
case is the most multi-dimensional.  
This case does not correspond to K3 surfaces but it is   
important. Considering a primitive sublattice $T_1\subset T$ with 
two positive squares and restriction of Borcherds 
form to $\Omega (T_1)$, one 
can construct other examples. But one should be very careful because 
this restriction may have additional zeros 
to quadratic divisors orthogonal to elements 
with square $-2$ of $T_1$ (e. g. it might be identically $0$).  
Considering this restriction, 
R. Borcherds found the form $\Phi$ with \thetag{2} for 
$T=U(2)\oplus U \oplus E_8(-2)$ (we denote by $K(t)$ a lattice $K$ with 
the form multiplied by $t\in \bq$).  
This case corresponds to moduli of K3 surfaces 
which cover twice Enriques surfaces. 
Considering orthogonal complement to $c\in U(2)$, 
we get the mirror symmetric family 
with $S=U\oplus E_8(-2)$. This corresponds to K3 surfaces with 
involution having the set of fixed points 
equals to union of two elliptic curves 
(see \cite{N3} or \cite{N8}). Both these cases are parabolic.   

In our papers \cite{GN1}---\cite{GN6} we mainly considered 
the case when $\rk S=3$. Respectively, $\rk T=5$. In particular, for 
$$
T=2U\oplus \langle -2t \rangle,\ \ t =1,\, 2,\, 3,\, 4,
\tag{3}
$$
and 
$$
T=2U(k)\oplus \langle -2 \rangle, \ \ k=1,\, \dots,\, 8,\, 10,\, 12,\, 16,  
\tag{4} 
$$
we found $2$-reflective automorphic forms $\Phi$ with expansion  
of type \thetag{2} for an isotropic $c$ 
in the first summand $U$ or $U(k)$ respectively. Then  
$$
S=c^\perp_T/\bz c =
\cases
U\oplus \langle -2t \rangle &\text{for $T=2U\oplus \langle -2t \rangle$,}\\
U(k)\oplus \langle -2 \rangle &\text{for $T=2U(k) \oplus \langle -2 \rangle$.}
\endcases
\tag{5} 
$$
Thus, the families $\M_S$ and $\M_{T^\perp}$ are mirror symmetric 
families of K3 surfaces for our mirror symmetry. 
It seems, that for many lattices $S$ and $T$ from \thetag{5}  
existence of Calabi--Yau manifolds with K3 fibers having these Picard and 
transcendental lattices is not known. Problems A, B and 
A$^\prime$, B$^\prime$  are very interesting 
for these $S$ and $T$.    

Physicists we have mentioned in Introduction   
considered several examples when $2$-reflective automorphic 
forms $\Phi$ take part in calculation of the Yukawa coupling (for model B) 
and the quantum intersection pairing (for model A) when 
existence of Calabi--Yau manifolds fibrated by the corresponding 
K3 surfaces was known.   

T. Kawai in \cite{Ka2} considered $S=U\oplus \langle -2 \rangle$ and 
$T=2U\oplus \langle -2 \rangle$.  
For the model (A) he took Calabi--Yau 3-folds of degree 
$20$ in the weighted projective 
space ${\Bbb P}(10,3,3,2,2)$. They are naturally 
fibrated by K3 surfaces with the Picard lattice $S$.   
There are two $2$-reflective automorphic forms $\Phi$ on $\Omega (T)$ 
with respect to $O^+(T)$. One of them is the classical Siegel modular 
form $\Delta_5$ which is the product of 
even theta-constants. Another one 
is the well-known Igusa modular form $\Delta_{35}$ which is the 
first Siegel modular form of odd weight. 
For both these forms we found infinite product 
expansions of type \thetag{2} 
in \cite{GN1}, \cite{GN2}, \cite{GN4}. T. Kawai used 
combination of these forms for the quantum intersection pairing of  
the Calabi--Yau 3-folds above.  
   
The Borcherds automorphic form for 
$T=U\oplus U(2)\oplus \langle -2 \rangle$  
(we discussed it above) was recently used by J. Harvey 
and G. Moore in \cite{HM3}. For model (B)    
they used Calabi--Yau 3-folds constructed by 
C. Borcea \cite{Bo} and Cl. Voisin \cite{V}. They are fibrated by 
K3 surfaces with the transcendental lattice $T$.

One can suggest that: {\it The $2$-reflective 
automorphic forms of type \thetag{2} used for 
the variant of mirror symmetry for K3 surfaces described above 
always take part in the quantum intersection pairing or Yukawa 
coupling of Calabi--Yau manifolds fibrated by the corresponding K3 
surfaces if the Calabi--Yau manifolds do exist.} 
Certainly, one need to make this conjecture much more concrete. 

\vskip10pt 

\head 
4. Example of one of the most remarkable  
$2$-reflective automorphic forms in dimension $3$ 
\endhead 

We finish with an example of an automorphic form of 
type \thetag{2} from \cite{GN6}. We give it for the lattice 
$T=2U(12)\oplus \langle -2 \rangle$ (it is one of the lattices 
\thetag{4}). 

Using methods from \cite{G1}---\cite{G5}, 
we constructed in \cite{GN6} an  automorphic cusp form   
$\Delta_1$ of the minimal possible weight $1$ with respect 
to the orthogonal group of $T=2U(12)\oplus \langle -2\rangle$. 
It has a character of order $6$. We use basis 
$f_2, \hat{f}_3, f_{-2}$ of the lattice 
$S=U(12)\oplus \langle -2 \rangle$ with the Gram matrix 
$$
\left(\matrix
0&0&12\\
0&-2&0\\
12&0&0
\endmatrix
\right)
$$
and corresponding coordinates 
$z_1, z_2, z_3$. Then  
$$
\align
\Delta_1(z_1, z_2, z_3)&=
\sum_{M\ge 1}
\sum
\Sb
m >0,\,l\in \bz\\
\vspace{0.5\jot} n,\,m\equiv 1\,mod\,6\\
\vspace{0.5\jot}
4nm-3l^2=M^2
\endSb
\hskip-4pt
\biggl(\dsize\frac{-4}{l}\biggr)
\biggl(\dsize\frac{12}{M}\biggr)
\sum\Sb a|(n,l,m)\endSb \biggl(\dsize\frac{6}{a}\biggr)
q^{n/6}r^{l/2}s^{m/6}\\
{}&=
q^{1/6}r^{1/2}s^{1/6}
\prod
\Sb n,\,l,\,m\in \Bbb Z\\
\vspace{0.5\jot}
(n,l,m)>0
\endSb
\bigl(1-q^{n} r^{l} s^{m}\bigr)^{f_{3}(nm,l)}
\tag{6}
\endalign 
$$
where 
$q=\exp{(24\pi i z_1)}$,
$r=\exp{(4\pi i z_2)}$,
$s=\exp{(24\pi i z_3)}$
and 
$$
\biggl(\frac{-4}{l}\biggr)=\cases \pm 1 &\text{if }
l\equiv \pm 1\ \hbox{mod}\ 4\\
\hphantom{\pm}0 &\text{if }
l\equiv \hphantom{\pm} 0 \ \hbox{mod}\ 2
\endcases,\ \  
\biggl(\frac{12}M\biggl)=
\cases
\hphantom{-}1\  \text{ if }\  M\equiv \pm 1 \operatorname{mod} 12\\
-1\ \text{ if }\  M\equiv \pm 5 \operatorname{mod} 12\\
\hphantom{-}0\  \text{ if }\  (M,12)\ne 1
\endcases, 
$$
$$ 
\biggl(\frac{\,6\,}{\,a\,}\biggr)=\cases \pm 1 &\text{if }
a\equiv \pm 1\ \hbox{mod}\ 6\\
\hphantom{\pm}0 &\text{if }
(a,6)\not=1\ . 
\endcases 
$$
The multiplicities $f_3(nm,l)$ of the infinite product are 
defined by a weak Jacobi form
$\phi_{0,3}(\tau,z)=\dsize\sum_{n\ge 0,\,l\in \bz}f_3(n,l)q^nr^l$ 
of weight $0$ and index $3$ with integral Fourier coefficients: 
$$
\phi_{0,3}(\tau ,z)=
\biggl(\frac {\vartheta(\tau ,2z)}{\vartheta(\tau ,z)}\biggr)^2
=r^{-1}
\biggl(\prod_{n\ge 1}(1+q^{n-1}r)(1+q^{n}r^{-1})(1-q^{2n-1}r^2)
(1-q^{2n-1}r^{-2})\biggr)^2 
$$    
where $q=\exp{(2\pi i \tau )}$, $r=\exp{(2\pi i z)}$. 
The divisor of $\Delta_1$ is sum with multiplicities one 
of all quadratic divisors orthogonal to elements of 
$T$ with square $-2$. The $\Delta_1$ 
defines the generalized Lorentzian Kac--Moody superalgebra 
with the denominator function \thetag{6}  
(see \cite{GN3} and \cite{GN1}---\cite{GN6} 
for details and other examples). 
Conjecturally the algebra is related with symmetries 
of some physical theory.  

It seems that for the automorphic form $\Delta_1$  
existence of Calabi--Yau manifolds fibrated 
by K3 surfaces with the corresponding Picard lattice 
$S=U(12)\oplus \langle -2 \rangle$ or the transcendental 
lattice $T=2U(12)\oplus \langle -2 \rangle$ are not known.

\enddocument
\end